\documentclass[12pt,a4paper]{amsart}

\newcommand{\N}{\rm{I\!N}}

\newtheorem{theorem}{Theorem}

\newtheorem{prop}[theorem]{Proposition}

\newcommand{\eps}{\varepsilon}

\newcommand{\ohne}{\setminus}
\newcommand{\gen}{\rightarrow}
\newcommand{\aus}{\subset}
\newcommand{\Norm}[1]{\Bigl\|#1\Bigr\|}
\newcommand{\norm}[1]{\|#1\|}

\newcommand{\betr}[1]{| #1  |}
\newcommand{\eing}[1]{_{|{#1}}}
\newcommand{\ebew}{\hfill{\rule{1.2ex}{1.2ex}}}
\newcommand{\bgl}{\begin{eqnarray}}
\newcommand{\bglst}{\begin{eqnarray*}}
\newcommand{\egl}{\end{eqnarray}}
\newcommand{\eglst}{\end{eqnarray*}}
\newcommand{\Pel}{Pe\l\-czy\'ns\-ki}
\newcommand{\Ref}[1]{(\ref{#1})}

\newcommand{\wst}{$w^{*}$}
\newcommand{\mdE}{|\;}

\begin{document}
\title{The dual of a non-reflexive L-embedded Banach space contains $l^{\infty}$ isometrically}
\author{Hermann Pfitzner}
\address{Universit\'e d'Orl\'eans\\
BP 6759\\
F-45067 Orl\'eans Cedex 2\\France}
\email{hermann.pfitzner@univ-orleans.fr}
\begin{abstract}
See title. (A Banach space is said to be L-embedded if it is complemented in its bidual such that the norm between
the two complementary subspaces is additive.)
\end{abstract}
\maketitle\noindent
This note is an afterthought to a result of Dowling \cite{Dow-asympt} according to which a dual Banach space contains
an isometric copy of $c_0$ if it contains an asymptotic one. (For definitions see below.)
It is known (\cite{Pf-L1} or \cite[Th. IV.2.7]{HWW}) that the dual of a non-reflexive L-embedded Banach  space contains $c_0$ isomorphically.
For a special class of L-embedded Banach  spaces the construction of the $c_0$-copy has been
improved so to yield an asymptotic one (\cite[Prop.\ 6]{Pfi-Fixp}) and it turns out that this improvement is possible in the general case
which together with Dowling's result yields isometric copies of $c_0$ in the dual of an L-embedded Banach  space.

All this is perhaps known - or at least not surprising - to experts in the field but the final result,
i.e. an isometric copy of $c_0$ (and hence of $l^{\infty}$) in the dual of an L-embedded space $X$,
is optimal in the category of Banach spaces and therefore it seems worthwile proving it explicitely.
As in \cite{Pf-L1} we will prove a bit more by constructing the $c_0$-copy within the context of
\Pel's property (V$^*$) i.e. the $c_0$-basis will be constructed so to behave approximately like biorthogonal functionals
on the basis of a given $l^1$-basis in $X$, see \Ref{gl1th} and \Ref{gl2th} below where in particular the value $\tilde{c}_J(x_n)$ in \Ref{gl1th} is optimal.
(For the definition and some basic results on \Pel's property (V$^*$) see \cite{HWW}.)\medskip

\noindent{\em Preliminaries}:
A projection $P$ on a Banach space $Z$ is called an  {\em L-projection} if
$\norm{Pz}+\norm{z-Pz}=\norm{z}$ for all $z\in Z$.
A Banach space $X$ is called {\em L-embedded} (or {\em an L-summand in its bidual}) if it is the image of an L-projection on its bidual.
In this case we write $X^{**}=X\oplus_1 X_{\rm s}$.
Among classical Banach spaces, the Hardy space $H_0^1$, $L^1$-spaces and, more generally, the preduals of von Neumann algebras or of JBW$^*$-triples
serve as examples of L-embedded spaces.
A sequence $(x_n)$ in a Banach space $X$ is said to {\em span $c_0$
asymptotically isometrically} (or just {\em to span $c_0$ asymptotically})
if there is a null sequence $(\delta_n)$ in $[0,1[$ such that
$\sup(1-\delta_n)\betr{\alpha_n}
\leq \Norm{\sum\alpha_n {x_n}}
\leq \sup(1+\delta_n)\betr{\alpha_n}$
for all $(\alpha_n)\in c_0$.
$X$ is said to contain $c_0$ asymptotically if it contains such a sequence $(x_n)$.
Recall the routine fact that if $(x_n^*)$ in $X^*$ is equivalent to the canonical
basis of $c_0$ then $\sum\alpha_n x_n^*$ makes
sense for all $(\alpha_n)\in l^{\infty}$ in the \wst-topology of $X^*$
and by lower \wst-semicontinuity of the norm an estimate $\norm{\sum\alpha_n x_n^*}\leq M\sup\betr{\alpha_n}$ that holds for
all $(\alpha_n)\in c_0$ extends to all $(\alpha_n)\in l^{\infty}$.
The Banach spaces we consider in this note are real or complex, the set $\N$ starts at $1$.

To a bounded sequence $(x_{n})$ in a Banach space $X$ we associate its 'James constant'
\bglst
c_J (x_{n}) =\sup c_m \quad\mbox{ where the }\quad
c_m=\inf_{\sum_{n\geq m}\betr{\alpha_n} =1}\norm{\sum_{n\geq m} \alpha_n x_n}
\eglst
form an increasing sequence.
If $(x_{n})$ is equivalent to the canonical basis of $l^1$ then $c_J(x_{n})>0$
and more specifically, $c_J(x_{n})>0$ if and only if there is an integer $m$ such that $(x_{n})_{n\geq m}$ is
equivalent to the canonical basis of $l^1$.
Roughly speaking, the number $c_J (x_{n})$ may be thought of as the 'approximately best $l^1$-basis constant' of $(x_{n})$;
precisely speaking, there is a null sequence $(\tau_m)$ in $[0,1[$ (determined by $c_m=(1-\tau_m)c_J(x_n)$) such that
$\norm{\sum_{n=m}^{\infty} \alpha_n x_n}\geq (1-\tau_m)c_J(x_n)\sum_{n=m}^{\infty}\betr{\alpha_n}$ for all $(\alpha_n)\in l^1$
and $c_J (x_{n})$ cannot be replaced by a strictly greater constant.
If one passes to a subsequence $(x_{n_k})$ of $(x_n)$ then $c_J(x_{n_k})\geq c_J(x_n)$ hence it makes sense to define
$$\tilde{c}_J(x_n)=\sup_{n_k}c_J(x_{n_k}).$$
The standard reference for L-embedded Banach spaces is the monograph \cite[Ch.\ IV]{HWW}.
For general Banach space theory and undefined notation we refer to \cite{Die-Seq}, \cite{JohLin}, or \cite{LiTz12}.\bigskip

\noindent
The main result of this note is
\begin{theorem}
Let $X$ be an L-embedded Banach space and let $(x_n)$ be equivalent to the canonical basis of $l^1$.
Then there is a sequence $(x_n^*)$ in $X^*$ that generates $l^{\infty}$ isometrically, more precisely
\bgl
\norm{\sum\alpha_n x_n^*}= \sup\betr{\alpha_n}\mbox{ for all } (\alpha_n)\in l^{\infty}\label{gl3prop}
\egl
and there is a strictly increasing sequence $(p_n)$ in $\N$ such that
\bgl
\lim \betr{x_n^*(x_{p_n})} &=& \tilde{c}_J(x_m)\label{gl1th}\\
x_n^*(x_{p_l}) &=&0 \quad\quad\mbox{ if } l<n \label{gl2th}
\egl
In particular, the dual of a non-reflexive L-embedded Banach space contains an isometric copy of $l^{\infty}$.
\end{theorem}
\noindent
In order to prove the theorem we first state and prove Dowling's result in a way which fits our purpose.
\begin{prop}
Let $(\eps_n)$ be a null sequence in $[0,1[$,
let $(N_n)$ be a sequence of pairwise disjoint infinite subsets of $\N$
 and let
$(y_n^*)$ in the dual of a Banach space $Y$ span $c_0$ such that
\bgl
\norm{\sum\alpha_n y_n^*}\leq \sup(1+\eps_n)\betr{\alpha_n}
\quad\mbox{ and }\quad
\norm{y_n^*}\gen 1\label{gl1prop}
\egl
for all $(\alpha_n)\in c_0$.

Then the elements
\bgl
x_n^*=\sum_{k\in N_n}\frac{y_k^*}{1+\eps_k}\label{gl2prop}
\egl
generate $l^{\infty}$ isometrically (as in \Ref{gl3prop}).
\end{prop}
\noindent{\em Proof of the proposition}:
With $(\eps_n)$, $(N_n)$ and $(y_n^*)$ as in the hypothesis of the statement define $x_n^*$ by \Ref{gl2prop}.
Then $\norm{x_n^*}\leq1$ for all $n\in\N$ by the first half of \Ref{gl1prop}. For the inverse inequality we have that
\bglst
\norm{x_n^*}\geq\norm{2\frac{y_m^*}{1+\eps_m}}-\Norm{\frac{y_m^*}{1+\eps_m}-\sum_{k\in N_n, k\neq m}\frac{y_k^*}{1+\eps_k}}
\geq 2\frac{\norm{y_m^*}}{1+\eps_m}-1
\eglst
holds for all $m\in N_n$ hence $\norm{x_n^*}\geq1$ by the second half of \Ref{gl1prop} which proves $\norm{x_n^*}=1$.

Similarly we show \Ref{gl3prop}: First, ``$\leq$'' of  \Ref{gl3prop} follows from the first half of \Ref{gl1prop};
second, by the just shown inequality we have
\bglst
\norm{\sum\alpha_n x_n^*}
\geq 2\betr{\alpha_m} - \norm{\alpha_m x_m^*-\sum_{n\neq m}\alpha_n x_n^*}
\geq 2\betr{\alpha_m} -\sup\betr{\alpha_n}
\eglst
for all $m\in\N$ hence ``$\geq$'' of  \Ref{gl3prop}.
\ebew\bigskip\\
\noindent{\em Proof of the theorem}:\\
Let $(\delta_n)$ be a sequence in $]0,1[$ converging to $0$. Suppose $(x_n)$ is an $l^1$-basis and write $\tilde{c}=\tilde{c}_J(x_n)$ for short.\\
{\em Observation:}
Given $\tau>0$ there is a subsequence $(x_{n_k})$ of $(x_n)$ such that $c_J(x_{n_k})>(1-\tau)\tilde{c}$ and
by James' $l^1$-distortion theorem there are blocks of the $x_{n_k}$ which span $l^1$ almost isometrically that is to say
there are pairwise disjoint finite sets $A_l\subset\{n_k\mdE k\in\N\}$, a sequence of scalars $(\lambda_n)$
such that $\sum_{k\in A_l}\betr{\lambda_k}=1$
and such that the sequence $(z_l)$ defined by $z_l=\tilde{z}_l/\norm{\tilde{z}_l}$ and  $\tilde{z}_l=\sum_{k\in A_l}\lambda_k x_k$
satisfies
$(1-2^{-m})\sum_{l=m}^{\infty}\betr{\alpha_l}\leq
\Norm{\sum_{l=m}^{\infty} \alpha_l {z_l}}   \leq (1+2^{-m})\sum_{l=m}^{\infty}\betr{\alpha_l}$ for all $m\in\N$;
furthermore $\norm{\tilde{z}_l}\gen c_J(x_{n_k})$ whence the existence of $l'$ such that $\betr{\tilde{c}-\norm{\tilde{z}_{l'}}}<\tau$.

By induction over $n\in\N$ we will construct finite sequences $(y_{i}^{(n)*})_{i=1}^n$ in $X^*$,
a sequence $(\tilde{y}_n)$ in $X$, pairwise disjoint finite sets $C_n\subset\N$ and a scalar sequence $(\mu_n)$
such that, with the notation $y_n=\tilde{y}_n/\norm{\tilde{y}_n}$,
\bgl
\sum_{k\in C_n}\betr{\mu_k}=1,&&  \tilde{y}_n=\sum_{k\in C_n}\mu_k x_k,\quad \betr{\tilde{c}-\norm{\tilde{y}_n}}<\delta_n,
\label{gl6a}\\
\betr{y_{i}^{(n)*}(y_i)}
&>& 1-\delta_i  \quad\quad\quad\quad\quad\quad\quad\quad\quad\forall i\leq n,
\label{gllm5}\\
y_{i}^{(n)*}(y_l)
&=&0  \quad\quad\quad\quad\quad\quad\quad\quad\quad\quad\quad\forall l<i\leq n,
\label{gl1}\\
y_{i}^{(n)*}(x_p)
&=&0  \quad\quad\quad\quad\quad\quad\quad\quad\quad\quad\quad\forall p\in C_l,\;\forall l<i\leq n,
\label{gl1bis}\\
\Norm{\sum_{i=1}^{m}\alpha_{i}y_{i}^{(n)*}}
&\leq&
\max_{i\leq m}(1+(1-2^{-n})\,\delta_i)\,\betr{\alpha_{i}}
\quad\forall m\leq n,\, \alpha_i \mbox{ scalars. }
\label{gllm6}
\egl
For $n=1$ we use the observation with $\tau=\delta_1$ and choose $l_1$ such that $\betr{\norm{\tilde{z}_{l_1}}-\tilde{c}}<\delta_1$.
Then we choose $y_{1}^{(1)*}$ such that $\norm{y_{1}^{(1)*}}=1$ and $y_{1}^{(1)*}(z_{l_{1}})=\norm{z_{l_{1}}}$.
It remains to set $C_1=A_{l_1}$, $\mu_k=\lambda_k$ for $k\in C_1$ and $\tilde{y}_1=\tilde{z}_{l_1}$.

For the induction step $n\mapsto n+1$ we recall that $(P^*)\eing{X^*}$
is an isometric isomorphism from $X^*$ onto $X_{{\rm s}}^{\bot}$, that
$X^{***}=X^{\bot}\oplus_{\infty} X_{{\rm s}}^{\bot}$ and that 
$(P^*x^*)\eing{X}=(x^*)\eing{X}$ for all $x^*\in X^*$.
Let $(z_l)$ be as in the observation above with $\tau=\delta_{n+1}$ and
let $z_{{\rm s}}\in X^{**}\ohne X$ be a \wst-accumulation point of the $z_l$.
Then $z_{{\rm s}} \in X_{{\rm s}}$ and $\norm{z_{{\rm s}}}=1$ by the proof of \cite[Lem\ 1]{Pfi-Fixp} (or by some general folklore argument).
Choose $t\in {\rm ker}\, P^*\aus X^{***}$ such that $\norm{t}=1$ and $t(z_{{\rm s}})=\norm{z_{{\rm s}}}$.
Put
\bglst
E&=&{\rm lin}(\{P^*y_{i}^{(m)*}\mdE  i\leq m\leq n\}\cup \{t\})\aus X^{***},\\
F&=&{\rm lin}(\{y_i\mdE i\leq n\}\cup\{z_{{\rm s}}\}\cup\{x_p\mdE p\in\bigcup_{l\leq n}C_l\})\aus X^{**}
\eglst
and choose $\eta>0$ such that
\bglst
(1+\eta)  (1+(1-2^{-n})\,\delta_i) &<& 1+(1-2^{-(n+1)})\,\delta_i \mbox{ and }\eta < (1-2^{-(n+1)})\,\delta_{n+1}
\eglst
for all $i\leq n$.
The principle of local reflexivity provides an operator $R:E\gen X^*$ such that
\bgl
(1-\eta)\norm{e^{***}}&\leq&\norm{Re^{***}}
\leq(1+\eta)\norm{e^{***}},\label{gl7a}\\
f^{**}(Re^{***})&=&e^{***}(f^{**}),\label{gl8a}
\egl
for all $e^{***}\in E$ and $f^{**}\in F$.\\
We define $y_{i}^{(n+1)*}=R(P^*y_{i}^{(n)*})$ for $i\leq n$ and
$y_{n+1}^{(n+1)*}=Rt$ and obtain (\ref{gllm6}, $n+1$) (with $\alpha_i=0$ if $m<i\leq n+1$) by
\bglst
\Norm{\sum_{i=1}^{n+1}\alpha_{i}y_{i}^{(n+1)*}}
&\stackrel{\Ref{gl7a}}{\leq}&
(1+\eta)\,\Norm{\Bigl(\sum_{i=1}^{n}\alpha_{i}P^*y_{i}^{(n)*}\Bigr)+
\alpha_{n+1}t}\\
&=&
(1+\eta)\max\Bigl(\Norm{
\sum_{i=1}^{n}\alpha_{i}P^*y_{i}^{(n)*}},
\norm{\alpha_{n+1}t}\Bigr)\\
&=&
(1+\eta)\max\Bigl(\Norm{
\sum_{i=1}^{n}\alpha_{i}y_{i}^{(n)*}},
\norm{\alpha_{n+1}t}\Bigr)\\
&\stackrel{(\ref{gllm6})}{\leq}&
(1+\eta)\max\Bigl(\max_{i\leq n}(1+(1-2^{-n})\,\delta_i)\betr{\alpha_{i}},\betr{\alpha_{n+1}}\Bigr)\\
&\leq &
\max_{i\leq n+1}(1+(1-2^{-(n+1)})\,\delta_i)\betr{\alpha_{i}}.
\eglst
Since $z_{{\rm s}}$ is a \wst-cluster point of $(z_l)$ we have
\bglst
\betr{y_{n+1}^{(n+1)*}(z_l)}
&>&
\betr{z_{{\rm s}}(y_{n+1}^{(n+1)*})}
-\delta_{n+1}\\
&\stackrel{\Ref{gl8a}}{=}&
\betr{t(z_{{\rm s}})}-\delta_{n+1}=1-\delta_{n+1}
\eglst
for infinitely many $l$;
furthermore, an $l_{n+1}$ can be chosen among these $l$ so to obtain $\betr{\norm{\tilde{z}_{l_{n+1}}}-\tilde{c}}<\delta_{n+1}$.
Set $C_{n+1}=A_{l_{n+1}}$, $\tilde{y}_{n+1}=\tilde{z}_{l_{n+1}}$, $\mu_k=\lambda_k$ for $k\in C_{n+1}$.
Then \Ref{gl6a} holds and (\ref{gllm5}, $n+1$) holds for $i=n+1$.
For $i\leq n$, (\ref{gllm5}, $n+1$) follows from
$$y_{i}^{(n+1)*}(y_i)=(P^*y_{i}^{(n)*})(y_i)=
y_{i}^{(n)*}(y_i)\stackrel{\Ref{gllm5}}{>}1-\delta_i.$$
Condition (\ref{gl1}, $n+1$) holds for $i=n+1$ by
\bglst
y_{n+1}^{(n+1)*}(y_l)=(Rt)(y_l)\stackrel{\Ref{gl8a}}{=}t(y_l)=0\quad\forall\, l<n+1
\eglst
and it holds for $i<n+1$ by
$$y_{i}^{(n+1)*}(y_l)=(P^*y_{i}^{(n)*})(y_l)=y_{i}^{(n)*}(y_l)\stackrel{\Ref{gl1}}{=}0\quad\forall\, l<i.$$
The proof of (\ref{gl1bis}, $n+1$) works like the one of (\ref{gl1}, $n+1$).
This ends the induction.

Now we define $y^*_i =\frac{1}{1+\delta_i}\lim_{n\in{\mathcal U}}y_{i}^{(n)*}$
for all $i\in\N$ where ${\mathcal U}$ is a fixed nontrivial ultrafilter on $\N$ and
where the limit is understood in the \wst-topology of $X^*$.
Then by \wst-lower semicontinuity of the norm and by \Ref{gllm6}
\bglst
\norm{\sum\alpha_{i}y^*_i}
\leq
\sup(1+\delta_i)\frac{\betr{\alpha_{i}}}{1+\delta_i}
=
\sup\betr{\alpha_{i}}
\eglst
for all $(\alpha_i)\in l^{\infty}$.
In particular, $\norm{y^*_i}\leq1$ hence $\norm{y_i^*}\gen1$ by \Ref{gllm5}
and $(y_i^*)$ satisfies \Ref{gl1prop} for $\eps_n=0$.

Let $(N_n)$ be a sequence of pairwise disjoint infinite subsets of $\N$
such that $(i_n)$ increases strictly where $i_n=\min N_n$.
By the proposition the sequence defined by
$$x_n^*=\sum_{i\in N_n}y_i^*$$
generates $l^{\infty}$ isometrically 
and we have $\betr{x_n^*(y_{i_n})}\stackrel{\Ref{gl1}}{=}\betr{y_{i_n}^*(y_{i_n})}\stackrel{\Ref{gllm5}}{\geq}1-\delta_{i_n}$.
By construction of the $y_i$ there is, for each $n\in\N$, an index $p_n\in C_{i_n}$ such that
$\betr{x_n^*(x_{p_n})}\geq (1-\delta_{i_n})\norm{\tilde{y}_{i_n}} \stackrel{\Ref{gl6a}}{\geq}(1-\delta_{i_n})(\tilde{c}-\delta_{i_n})$
which will yield ``$\geq$'' of \Ref{gl1th}.
In order to show ``$\leq$'' of \Ref{gl1th} suppose to the contrary that $x_{n_m}^*(x_{p_{n_m}})>\kappa+\tilde{c}$ for appropriate subsequences, all $m$ and
$\kappa>0$. According to an extraction lemma of Simons  \cite{Sim-DPP} we may furthermore suppose that
$\sum_{j\neq m} \betr{x_{n_j}^*(x_{p_{n_m}})}  <\kappa/2$ for all $m$. Then given $(\alpha_m)$ and $\theta_m$ such that
$\theta_m \alpha_m=\betr{\alpha_m}$ we obtain
$\norm{\sum\alpha_m x_{p_{n_m}}}\geq(\sum_j \theta_j x_{n_j}^*)(\sum_m \alpha_m x_{p_{n_m}})
\geq (\kappa+\tilde{c})\sum_m \betr{\alpha_m} - \sum_m \sum_{j\neq m}\betr{\alpha_m}\;\betr{x_{n_j}^*(x_{p_{n_m}})}
\geq((\kappa/2) +\tilde{c})\sum_m \betr{\alpha_m}$
which yields the contradiction $c_J(x_{p_{n_m}})>\tilde{c}$
and thus shows ``$\leq$'' and all of \Ref{gl1th}
whereas \Ref{gl2th} follows from \Ref{gl1bis} via $y_i^*(x_p)=0$ for $p\in C_l$, $l<i$.

The last assertion of the theorem is immediate from the fact that non-reflexive L-embedded spaces contain $l^1$ isomorphically
\cite[IV.2.3]{HWW}\ebew\medskip

\noindent
Remarks:\\
1. It is not clear whether \Ref{gl2th} can be obtained also for $l>n$.
What can be said by Simons' extraction lemma (used in the proof) is that, under the assumptions of the theorem and given $\eps>0$,
it is possible (after passing to appropriate subsequences) to obtain in addition to \Ref{gl2th}
that $\sum_{n=1}^{l-1} \betr{x_n^*(x_{p_l})} =\sum_{n\neq l} \betr{x_n^*(x_{p_l})} <\eps$ for all $l$.
In case $\tilde{c}_J(x_n)=1=\lim \norm{x_n}$ (which happens when the $x_n$ span $l^1$ almost isometrically) this can be improved to
\bgl
\sum_{n\neq l} \betr{x_n^*(x_{p_l})}=  (\sum_n \betr{x_n^*(x_{p_l})})-\betr{x_l^*(x_{p_l})}\leq \norm{x_{p_l}}-\betr{x_l^*(x_{p_l})}\gen0.\label{glBem1}
\egl
One might also construct straightforward perturbations of the $x_n^*$ in order to get \Ref{gl2th} for $l\neq n$ but then it is not clear
whether these perturbations can be arranged to span $c_0$ isometrically, not just almost isometrically.

Since in general L-embedded spaces do not contain $l^1$ isometrically (see below, last remark) it is in general not possible to
improve \Ref{gl1th} and \Ref{gl2th} so to
obtain $x_n^*(x_{p_l})=\tilde{c}(x_m)$ if $l=n$ and $=0$ if $l\neq n$.\\
2.
As already alluded to in the  introduction, the construction of $c_0$ in this paper bears much resemblance to the one of \cite{Pf-L1}.
A different way to construct $c_0$ is contained in \cite{Pfi-unique} but it seems unlikely that this construction can be improved to
yield an isometric $c_0$-copy.\\
3. It follows from (\ref{gl1th}) that in L-embedded spaces the $\sup$ in the definition of $\tilde{c}_J$ is attained by the James constant
of an appropriate subsequence. For general Banach spaces this is not known although it can be shown by a routine diagonal argument that
each bounded sequence $(x_n)$ admits
a $c_J$-stable subsequence $(x_{n_k})$ (meaning that $\tilde{c}_J(x_{n_k})=c_J(x_{n_k})$) whose James constant is arbitrarily near to
$\tilde{c}(x_n)$.\\
4. Each normalized sequence $(x_n)$ in an L-embedded Banach space that spans $l^1$ almost isomorphically contains a subsequence each of whose
\wst-accumulation points in the  bidual attains its norm on the dual unit ball.
To see this let $(x_n^*)$ and $(x_{p_n})$ be the sequences given by the theorem and by Simons' extraction lemma (see \Ref{glBem1} above),
let $x_{{\rm s}}$ be a \wst-accumulation point of the $x_{p_n}$
and let $x^*=\sum x_n^*$; then $\norm{x^*}=1$ and on the one hand $\norm{x_{{\rm s}}}=1$ by \cite{Pfi-Fixp} and on the other hand
$x_{{\rm s}}(x^*)=\lim x^*(x_{p_n})\stackrel{\Ref{glBem1}}{=}\lim x_n^*(x_{p_n})\stackrel{\Ref{gl1th}}{=} c_J(x_n)=1$.

It would be interesting to know whether this remark holds for the whole sequence $(x_n)$ instead of only a subsequence $(x_{p_n})$.
A kind of converse follows from \cite[Rem. 2]{Pfi-unique} for separable $X$: If $x_{{\rm s}}\in X_{{\rm s}}$ attains its norm on the dual unit ball
then it does so on the sum of a wuC-series.\\
5.  Let us finally note that the presence of isometric $c_0$-copies in $X^*$ does not necessarily entail the presence of isometric copies
of $l^1$ in $X$ even if $X$ is the dual of an M-embedded Banach space.
This follows from  \cite[Cor.\ III.2.12]{HWW} which states that there is
an L-embedded Banach space which is the dual of an M-embedded space (to wit the dual of $c_0$ with an equivalent norm) which is strictly convex
and therefor{e does not
contain $l^1$ isometrically although it contains, as do all non-reflexive L-embedded spaces, $l^1$ asymptotically
(\cite{Pfi-Fixp}, see \cite{DJLT} for the definition of asymptotic copies and the difference to almost isometric ones).

\end{document}